\documentclass[12pt, a4paper]{amsart}
\usepackage{amsmath, amssymb, amsthm, amsfonts}
\usepackage{mathrsfs}
\usepackage{mathtools}
\usepackage[top=2.6cm, bottom=2.5cm, left=2.6cm, right=2.6cm]{geometry}

\theoremstyle{plain}
\newtheorem{thm}{Theorem}

\theoremstyle{remark}
\newtheorem{remark}{Remark}

\begin{document}
\title{On homotopy elements represented by quotients of Lie groups} 
\author{Haruo Minami}
\address{H. Minami: Professor Emeritus, Nara University of Education}
\email{hminami@camel.plala.or.jp}
\begin{abstract} 
Consider the quotient $G/B$ of a simple matrix Lie group $G$ by a subgroup 
$B$ isomorphic to a direct product of some of $S^1$s and $S^3$s such that 
its adjoint representation can be extended over $G$. Then it naturally inherits 
a stable framing from a twisted left invariant framing $\mathscr{L}^\alpha$ of $G$ where $\alpha$ is the realization of a complex representation of $G$. In this note we want to add some homotopy elements represented by such quotient framed manifolds to those presented in a table of [E. Ossa 1982].
\end{abstract}

\maketitle

\section{Introduction}

Let $G$ be a simply connected compact simple Lie group of dimension $d$. 
Regard $G$ as a framed manifold with the left invariant framing 
$\mathscr{L}$ and write $[G, \mathscr{L}]$ for its bordism class in $\pi_d^S$. 
Then we have a table in ~\cite{O}, which lists the homotopy elements represented by these bordism classes based on the results of ~\cite{A},  ~\cite{S},  ~\cite{W}, ~\cite{K} and ~\cite{BS}. These elements can also be arranged in the table of the first 21 stable homotopy groups of \!\!~\cite{T} and ~\cite{MT} (and also ~\cite{KM}) as follows.

\vspace{2mm}

\begin{itemize}
\item[] \ $\pi_1^S=\mathbb{Z}_2\eta$ \ \text{where} \ 
$\eta=[S^1, \mathscr{L}]$,
\item[] \ $\pi_2^S=\mathbb{Z}_2\eta^2$,
\item[] \ $\pi_3^S=\mathbb{Z}_{24}[SU(2), \mathscr{L}]$,
\item[] \ $\pi_4^S=0$,
\item[] \ $\pi_5^S=0$,
\item[] \ $\pi_6^S=\mathbb{Z}_2[SU(2), \mathscr{L}]^2$,
\item[] \ $\pi_7^S=\mathbb{Z}_{240} a[7]$,
\item[] \ $\pi_8^S=\mathbb{Z}_2[SU(3), \mathscr{L}]\oplus
\mathbb{Z}_2\eta a[7]$, 
\item[] \ $\pi_9^S=\mathbb{Z}_2\eta [SU(3), \mathscr{L}]\oplus
\mathbb{Z}_2\eta^2a[7]\oplus\mathbb{Z}_2a[9]$,
\item[] \ $\pi_{10}^S=\mathbb{Z}_3 [Sp(2), \mathscr{L}]\oplus
\mathbb{Z}_2\eta a[9]$, 
\item[] \ $\pi_{11}^S=\mathbb{Z}_{504}a[11]$,
\item[] \ $\pi_{12}^S=0$,
\item[] \ $\pi_{13}^S=\mathbb{Z}_3 [SU(2), \mathscr{L}]
[Sp(2), \mathscr{L}]$,
\item[] \ $\pi_{14}^S=\mathbb{Z}_2 a[7]^2\oplus\mathbb{Z}_2
[G_2, \mathscr{L}]$, 
\item[] \ $\pi_{15}^S=\mathbb{Z}_{480} a[15]\oplus
\mathbb{Z}_2[SU(4), \mathscr{L}]$ \ \text{where} \ $[SU(4), \mathscr{L}]=\eta [G_2, \mathscr{L}]$, 
\item[] \ $\pi_{16}^S=\mathbb{Z}_2a[16]\oplus\mathbb{Z}_2 \eta a[15]$,
\item[] \ $\pi_{17}^S=\mathbb{Z}_2 \eta a[16]\oplus\mathbb{Z}_2 
[SU(2), \mathscr{L}][G_2, \mathscr{L}]\oplus\mathbb{Z}_2\eta^2a[15]\oplus\mathbb{Z}_2 a[17]$,
\item[] \ $\pi_{18}^S=\mathbb{Z}_8 a[18]\oplus\mathbb{Z}_2 \eta a[17]$,
\item[] \ $\pi_{19}^S=\mathbb{Z}_{264} a[19]\oplus\mathbb{Z}_2 b[19]$,
\item[] \ $\pi_{20}^S=\mathbb{Z}_{24} a[20]$,
\item[] \ $\pi_{21}^S=\mathbb{Z}_2[Sp(3), \mathscr{L}]\oplus
\mathbb{Z}_2a[7]^3$.
\end{itemize}

In this note we want to add new homotopy elements represented with Lie groups to the above ones as many as possible by adopting elements constructed through the method of ~\cite{M}. Let $\lambda=\rho_\mathbf{R}$ be the realification of a faithful complex representation $\rho\colon G \to GL(n, \mathbb{C})$ and 
let $S\subset G$ denote a circle subgroup isomorphic to $U(1)\subset Sp(1)$ 
when identifying $Sp(1)=SU(2)$. Then we obtain the following. 

\begin{thm}  In the setting above we have
\begin{itemize} 
\item[\,\,(i)] \ $a[7] \ \,=[Sp(2)/Sp(1), (\mathscr{L}^\lambda)_{Sp(1)}]$ \ and 
\ $2a[7]=[SU(3)/S, \mathscr{L}_S]$, 
\item[\,(ii)] \ $a[9] \ \,=[Sp(2)/S, (\mathscr{L}^\lambda)_S]$ \ and \ 
$[Sp(2), \mathscr{L}^\lambda]=\eta a[9]$,
\item[(iii)] \ $a[11]=[SU(4)/SU(2)\times S, 
(\mathscr{L}^\lambda)_{SU(2)\times S}]$,
\item[\,(iv)] \ $a[15]=[Sp(3)/Sp(1)\times Sp(1), 
(\mathscr{L}^{2\lambda})_{Sp(1)\times Sp(1)}]$ \ and \  
$2a[15]=[SU(4), \mathscr{L}^\lambda]$,
\item[\,\,(v)] \ $a[19]=[Sp(3)/S\times S, \mathscr{L}_{S\times S}]$,
\item[\,(vi)] \ $a[18]=[Sp(3)/Sp(1), (\mathscr{L}^\lambda)_{Sp(1)}]$,
\item[(vii)] \ 
$a[17]=[Sp(3)/Sp(1)\times S, (\mathscr{L}^\lambda)_{Sp(1)\times S}]$ \ 
and \ $[Sp(3)/Sp(1), (\mathscr{L}^{2\lambda})_{Sp(1)}]=\eta a[17]$,
\item[(viii)] \ $a[20]=[Sp(3)/S, \mathscr{L}_S]$,
\item[\,(ix)] \ \,$b[19]=[Sp(3)/S\times S, (\mathscr{L}^\lambda)_{S\times S}]$,
\item[\,\,(x)] \,\,$a[16]=[Sp(3)/Sp(1)\times S\times S, 
(\mathscr{L}^\lambda)_{Sp(1)\times S\times S}]$.
\end{itemize}
\end{thm}

To begin with, we recall the construction of the quotient framed manifods $G/T^r$ with the framing $(\mathscr{L}^{k\lambda})_{T^r}$ in ~\cite{M} along with its extesions which appears in the equations of the theorem. Let $T^r=S\times\cdots\times S
\subset G$ ($r$ times) and as noted above, suppose that $\rho(S)$ is of the form 
$\{\mathrm{diag}(z, z^{-1})\!\mid\! z\in U(1)\}$ for each factor $S$ of $T^r$. Let 
$\rho$ be interpreted as a map modified such that when a colomn vector of 
$\rho(g)$ $(g\in G)$ is actually acted on by an element $\rho(z)$ 
$(z\in S)$, all its components $v$ including those belonging to other columns are converted to $|v|$ and then in parallel the negatives of those $v$ are converted to $- |v|$. If we write $\rho(g)_z$ for $\rho(g)$ with the thus modified columns, then we can view $\rho$ (resp. $\rho_\mathbb{R}$) as a map from $G/T^r$ to 
$GL(n, \mathbb{C})$ (resp. $GL(2n, \mathbb{R})$) by assigning $g\cdot z$ to 
$\rho(g)_z$. Below we gives a generalization of the construction of 
$(\mathscr{L}^{k\lambda})_{T^r}$ for $T^r$ with quaternionic circles as $S$.
.   

Let $B^r\subset G$ be a subgroup obtained by replacing $r_1$ $(0\le r_1\le r)$   circle components of $T^r$ by the unit quaternion group, written $Sp(1)$ or 
$SU(2)$. But in this case, in addition to the action of $\rho(z)$  
$(z\in U(1)\subset Sp(1))$, we need to take into account that of $\rho(zj)$  
$(zj\in Sp(1))$. These elements divide the action of 
$q=rz+ju\in Sp(1)$ ($r\ge 0$, $u\in\mathbb{C}$) into two parts. 
Hence observing the action of such elements we find that $\rho$ presents 
a map from $G/B^r$ to $GL(n, \mathbb{C})$ under the assignment 
$g\cdot q \to \rho(g)_z$ and so $\lambda$ can be viewed as 
a map from $G/B^r$ to $GL(2n, \mathbb{R})$.

By applying the argument for the circle case to the principal 
$B^r$-bundle $p\colon G\to G/B^r$ we have 
\begin{equation*} 
\underline{\mathbb{R}}^r\oplus t(G)\cong p^*(t(G/B^r))
\oplus\underline{\mathbb{H}}^{r_1}\oplus\underline{\mathbb{C}}^{r-r_1}
\end{equation*}
where $t(M)$ denotes the tangent bundle of $M$ and   
$\underline{\mathbb{K}}^s=G\times\mathbb{K}^s $ is the product bundle. 
Dividing this by the right action of $B^r$ we have  
\begin{equation*} 
\underline{\mathbb{R}}^r\oplus t(G)/B^r\cong t(G/B^r)
\oplus\underline{\mathbb{R}}^{2r+2r_1}.
\end{equation*}
Let $\underline{\mathbb{R}}^a\oplus t(G)\cong 
\underline{\mathbb{R}}^a\oplus (G\times t_e(G))$ be the isomorphism 
induced by $\mathscr{L}^{k\lambda}$ where $a=2kn-d$ and $t_e(G)$ is the tangent space at the identity $e\in G$. Then due to the fact that the right action of $B^r$ on $t_e(G)$ operates on the right-hand side as the adjoint action we have    
\[\underline{\mathbb{R}}^a\oplus t(G)/B^r\cong 
\underline{\mathbb{R}}^a\oplus (G\times_{Ad_{B^r}} t_e(G)).\]  
Hence, assuming that $Ad_{B^r}$ is stably extendable over $G$, 
from the above we have  
\[\underline{\mathbb{R}}^{r+d+s+s}\cong t(G/B^r)
\oplus\underline{\mathbb{R}}^{2r+2r_1+s}.\]
This represents the framing $(\mathscr{L}^{k\lambda})_{B^r}$. 
Here we write $(G/B^r, (\mathscr{L}^{k\lambda})_{B^r})$ for $G/B^r$ equipped with $(\mathscr{L}^{k\lambda})_{B^r}$ and $[G/B^r, (\mathscr{L}^{k\lambda})_{B^r}]$ for its framed bordism class as usual in the case when the extendability of $Ad_{B^r}$ over $G$ is guaranteed.

Let $G/B^r\subset \mathbb{R}^{d-r+N}$ be a smooth embedding for large $N$. If we denote by $\mathbb{R}^N$ the normal plane at $p(e)$, then for any $k\ge 0$ the framing $(\mathscr{L}^{k\lambda})_{B^r}$ induces an embedding 
$\tau_{k\lambda}\colon G/B^r\times \mathbb{R}^N\to \mathbb{R}^{d-r+N}$. 
By performing the Pontryagin–Thom construction of this embedding 
we obtain a based map
$\tau_{k\lambda}^*\colon S^{d-r+N}=(\mathbb{R}^{d-r+N})^+ \to 
(G/B^r\times \mathbb{R}^N)^+$ where $+$ denotes the one-point ompactification.
Write $\tau_{G/B^r, \,(\mathscr{L}^{k\lambda})_{B^r}}$ for the stable map of the composite  
\begin{equation*}
S^{d-r+N} \xrightarrow{\tau_{k\lambda}^*}
(G/B^r\times \mathbb{R}^N)^+\xrightarrow{{p'}^+}(\mathbb{R}^N)^+=S^N
\end{equation*}
where $p'\colon G/B^r\times \mathbb{R}^N\to \mathbb{R}^N$ is the canonical projection. Then we know that its homotopy class corresponds isomorphically to 
the framed bordism class of $(G/B^r, (\mathscr{L}^{k\lambda})_{B^r})$.

\begin{remark}
1) Let $\lambda_k$ be the $k$-th exterior power of the standard complex representations  of $SU(n)$ and $Sp(n)$. Then by ~\cite{Y} we know that the complexifications of their adjoint representations are given by  
\[c({Ad}_{SU(n)})=\lambda_1\lambda_{n-1}-1, \quad   
c({Ad}_{Sp(n)})=\lambda_1^2-\lambda_2.\]
2) Using these equations together with the relation $rc=2$ where $r$ is the 
realification we see that $Ad_{SU(2)}$ and $Ad_{Sp(1)}$ can be represented 
up to constant as the restrictions of real representaions of $SU(n)$ and $Sp(n)$, respectively, as follows.
\[Ad_{SU(n)}-(n-2)r(\lambda_1)\!\mid\! SU(2)=Ad_{SU(2)}, \  
Ad_{Sp(n)}-(n-1)r(\lambda_1)\!\mid\! Sp(1)=Ad_{Sp(1)}.\] 
3) From the above equations we find that the the problem of extendability of  $Ad_{SU(2)}$ and  $Ad_{Sp(1)}$ over $G/B^r$ can be solved by replacing $\mathscr{L}$ by $\mathscr{L}^{r_1\lambda}$ (cf. ~\cite{G}). 
\end{remark}

\section{Proof of Theorem}

Now we proceed with the proof of the theorem based on the above table.
\begin{proof}[Proof of $\mathrm{(i)}$]
Let 
\[Sp(2)/Sp(1)\to Sp(2)/Sp(1)\times Sp(1)=\mathbb{H}P^1\]
be the canonical principal $Sp(1)$-bundle where we identify 
$\mathbb{H}P^1=S^4$ and let $\zeta$ be the associated quaternionic line 
bundle. Then from Proposotion 5.1 of  ~\cite{LS} we have
\begin{equation*}
e'_\mathbb{R}([S(\zeta), \Phi_\zeta])=-1/240
\end{equation*}
where $\Phi_\zeta$ is the induced framing of the sphere bundle $S(\zeta)$. This allows us to take $a[7]=[S(\zeta), \Phi_\zeta]$. But by the definition of 
$(\mathscr{L}^\lambda)_{Sp(1)}$ and the construction of $\Phi_\zeta$ in ~\cite{LS} we see that these two framings are equivalent when viewed as 
$S(\zeta)=Sp(2)/Sp(1)$. Hence by replacing them we have 
\begin{equation*}
a[7]=[Sp(2)/Sp(1), (\mathscr{L}^\lambda)_{Sp(1)}].
\end{equation*}

On the other hand, by Theorem 1 of ~\cite{M} we also have
\begin{equation*} 
e'_\mathbb{R}([SU(3)/S, \mathscr{L}_S])=-1/120.
\end{equation*} 
since $e'_\mathbb{R}=e_\mathbb{C}$ in the present case because of  
$\dim SU(3)/S\equiv 7 \mod 8$ ~\cite{A}.
Comparing this with the above equation we see that 
\[2a[7]=[SU(3)/S, \mathscr{L}_S],\]
including the sign.
\end{proof}

\begin{proof}[Proof of $\mathrm{(ii)}$]
Suppose that $\tau_{Sp(2), \,\mathscr{L}^\lambda}\simeq c_\infty$ where $c_\infty$ is the constant map at the base point and consider the canonical principal fibration
\begin{equation*} 
Sp(1)\times Sp(1)\hookrightarrow Sp(2)\to Sp(2)/Sp(1)\times Sp(1)
=S^4.
\end{equation*}
Here we write $Sp(1)_1\times Sp(1)_2$ for the subgroup $Sp(1)\times Sp(1)\subset Sp(2)$ in order to distinguish the first $Sp(1)$ from the second one.  
If we assume that the framing $\mathscr{L}^\lambda$ of $Sp(2)$ behaves on 
$Sp(1)_1$ as $\mathscr{L}$, then it must operate on $Sp(1)_2$ as 
$\mathscr{L}^\lambda$. This allows us to replace the fiber $Sp(1)_1\times Sp(1)_2$ by $Sp(1)_1\times D_2$ in the fibration above where $D_2$ denotes the closed unit ball in $\mathbb{H}$ with $Sp(1)_2$ as its boundary, i.e. 
$\partial D_2=Sp(1)_2$. Let us write $E_{Sp(2)}\to S^4$ for the fibration thus obtained. Then the null homotopy of $\tau_{Sp(2), \,\mathscr{L}^\lambda}$ can be naturally extended over $E_{Sp(2)}$. Now since $\tau_{Sp(1)_2, \,\mathscr{L}}$ yields an element of order 24, it becomes in fact a null homotopy of 24 times 
$\tau_{Sp(2)/Sp(1)_1, (\mathscr{L}^\lambda)_{Sp(1)_1}}$.  This contradicts the result of (i) that $\tau_{Sp(2)/Sp(1)_1, \, (\mathscr{L}^\lambda)_{Sp(1)_1}}$ represents $a[7]$. Hence we have $\tau_{Sp(2), \,\mathscr{L}^\lambda}\not\simeq c_\infty$; i.e.  
\begin{equation*}
[Sp(2), \mathscr{L}^\lambda]\ne 0.
\end{equation*}

Let $\xi$ be the complex line bundle associated to 
the canonical principal circle bundle
\begin{equation*}
S\hookrightarrow Sp(2)\to Sp(2)/S. 
\end{equation*}
Let  $\pi\colon Sp(2)/S\to Sp(2)/Sp(1)$ be the canonical bundle projection where we identify $Sp(2)/Sp(1)=S(\mathbb{H}\oplus\mathbb{H})$, the unit sphere in 
$\mathbb{H}\oplus\mathbb{H}$. Then we see that the restriction of $\xi$ on $Sp(1)/S$ is stably tivial due to the action of $\lambda$. However, making a further observation of the behaivor of $\lambda$, we have that the stable triviality of $\xi$ can be extended over $\pi^{-1}(S(\{0\}\oplus\mathbb{H}))$ and therefore, 
as a result, that $\xi$ itself becomes stably trivial. Based on this, applying the defining equation of Lemma 2.2 in ~\cite{O} to the sphere bundle $S(\xi)$, we have 
$[Sp(2), \mathscr{L}^\lambda]=\eta [Sp(2)/S, (\mathscr{L}^\lambda)_S]$. This  
together with the non zeroness result of it above allows us to take 
$a[9]=[Sp(2)/S, (\mathscr{L}^\lambda)_S]$. 
\end{proof}

\begin{proof}[Proof of $\mathrm{(iii)}$]
Let $SU(2)$ and $S$ be the subgroups of $SU(4)$ consisting of the 
elements 
\begin{equation*}
\left(
\begin{array}{cc} 
A & 0 \\ 
0 & I_2
\end{array}
\right) \ \text{with} \ A\in SU(2) \quad \text{and} \quad
\left(
\begin{array}{ccc} 
I_2 & 0 & 0 \\ 
0  & z  & 0  \\
0  & 0  & \bar{z}
\end{array}
\right) \ \text{with} \ z\in U(1),
\end{equation*}
respectively, where $I_k$ is the unit matrix of order $k$, and put 
$H=SU(2)\times S$ in $SU(4)$. 
We also now proceed in a similar manner 
to the proof of Theorem 1 in ~\cite{M}. Let us put 
\begin{equation*}
R_1(r_1z_1, u_1)=
\left(
\begin{array}{cccc} 
r_1z_1& 0 & u_1 & 0 \\ 
 0 & 1 & 0 & 0 \\
-\bar{u}_1& 0 & r_1\bar{z}_1 & 0\\
 0 & 0 & 0 & 1
\end{array}
\right), \ \
R_2(r_2z_2, u_2)=
\left(
\begin{array}{cccc} 
r_2z_2& 0 & 0 & u_2 \\ 
 0 & 1 & 0 & 0 \\
0& 0 & 1 & 0\\
-\bar{u}_2 & 0 & 0 & r_2\bar{z}_2
\end{array}
\right)
\end{equation*}
\begin{equation*}
R_3(r_3z_3, u_3)=
\left(
\begin{array}{cccc} 
1& 0 & 0 & 0 \\ 
 0 & r_3z_3 & u_3& 0 \\
0& -\bar{u}_3 & r_3\bar{z}_3 & 0\\
 0 & 0 & 0 & 1
\end{array}
\right), \ \
R_4(r_4z_4, u_4)=
\left(
\begin{array}{cccc} 
1& 0 & 0 & 0 \\ 
 0 & r_4z_4 & 0 & u_4 \\
0& 0 & 1 & 0\\
0& -\bar{u}_4  & 0 & r_4\bar{z}_4
\end{array}
\right)
\end{equation*}
\begin{equation*}
R_5(r_5z_5, u_5)=
\left(
\begin{array}{cccc} 
1& 0 & 0 & 0 \\ 
 0 & 1 & 0& 0 \\
0& 0 & r_5z_5 & 0\\
 0 & 0 & -\bar{u}_5 & r_5\bar{z}_5
\end{array}
\right) \qquad  \text{with} \
\left(
\begin{array}{cc} 
r_iz_i & u_i \\
 -\bar{u}_i  & r_i\bar{z}_i
\end{array}
\right)\in SU(2) \ \qquad 
\end{equation*}
in $SU(4)$ where $r_i\ge 0$ and $z_i \in U(1)$.
For any $z\in U(1)$ we define $d(z)\in SU(4)$ to be the product element
$R_1(z, 0)R_2(z, 0)R_3(z, 0)R_4(z, 0)R_5(z, 0)$. Put $P=SU(4)/H$. Then 
$d(z)$ acts on $P$ freely. Letting $S'=\{d(z)\mid z\in U(1)\}$ consider the projection $\pi\colon P\to M$ to be the principal $S'$-bundle over $M=P/S'$. 
Now by definition we see that 
$\textstyle\bigr(\prod_{i=1}^5R_i(r_iz_i, u_i)\bigl)\bigr(\prod_{i=1}^5d(\bar{z}_i)\bigl)$ 
can be written in the form $\prod_{i=1}^5R_i(r_i, w_i)$. Therefore, writing 
$S^2=\{(r, w)\mid r\ge 0, w\in\mathbb{C}, r^2+|w|^2=1 \}$, we have a map 
\[\phi\colon S^2\times \cdots\times S^2  \to M\qquad \ \  
\mathrm{(five \ times)}\]
given by $((r_1, w_1), \ldots, (r_5, w_5)) \to \pi\bigr(R_1(r_1, w_1)\cdots 
R_5(r_5, w_5)\bigl)$. Let $E$ be the canonical complex line bundle associated to 
$\pi$. Then we see that $\phi^*E$ can be written as the exterior tensor product 
$L\boxtimes\cdots\boxtimes L$ of five the complex Hopf line bundle 
$L$ over $S^2$. By verifying that the restriction of $\phi$ to 
$(S^2-\{\infty\})\times \cdots\times (S^2-\{\infty\})$ is injective as is seen in ~\cite{M} we find that it can be deformed to a degree one map. Applying these two facts to Proposition 2.1 of ~\cite{LS} we obtain 
$e_\mathbb{C}([SU(4)/H, (\mathscr{L}^\lambda)_H])=-1/252$, so we have
\[e'_\mathbb{R}([SU(4)/H, (\mathscr{L}^\lambda)_H])=-1/504\] in this case since 
$\dim SU(4)/H\equiv 3\mod 8$. This shows that (iii) holds true.
\end{proof}

\begin{proof}[Proof of $\mathrm{(iv)}$]
Let 
\[Sp(1)\hookrightarrow Sp(3)/Sp(1)\times Sp(1)\to 
Sp(3)/Sp(1)\times Sp(1)\times Sp(1)\]
be the canonical principal $Sp(1)$-bundle and let $\zeta$ be the associated quaternionic line bundle. Then 
applying Proposition 5.1 of ~\cite{LS} we have $a[15]=[S(\zeta), \Phi_\zeta]$
as in the case of $a[7]$. Reinterpreting this in a similar way we obtain 
\begin{equation*}
a[15]=[Sp(3)/Sp(1)\times Sp(1), (\mathscr{L}^{2\lambda})_{Sp(1)\times Sp(1)}].
\end{equation*}

Moreover, by Theorem 2 of ~\cite{M} we also have 
$e_\mathbb{C}([SU(4), \mathscr{L}^\lambda])=-1/240$ and so  
\begin{equation*}
2a[15]=[SU(4), \mathscr{L}^\lambda]
\end{equation*}
for the same reason as in the case of (i).
\end{proof}

\begin{proof}[Proof of $\mathrm{(v)}$]
By Theorem 1 of ~\cite{M} we have
$e_\mathbb{C}([Sp(3)/T^2, \mathscr{L}_{T^2}])=5/132$ and so 
\[e'_\mathbb{R}([Sp(3)/T^2, \mathscr{L}_{T^2}])=5/264.\]
This means that $[Sp(3)/T^2, \mathscr{L}_{T^2}]$ satisfies 
the required condition for being $a[19]$. 
\end{proof}

\begin{proof}[Proof of $\mathrm{(vi)}$]
Let us consider the canonical principal fibration  
\[Sp(2)\hookrightarrow Sp(3)/Sp(1)\to Sp(3)/Sp(1)\times Sp(2)
=\mathbb{H}P^2.\]
Here $\mathbb{H}P^2=S^4\cup_\nu D^8$ where $D^8$ denotes the hemisphere consisting of the first column elements $(q_1, q_2, q_3)$ of $Sp(3)$ 
with $q_3=r\ge 0$ and  $\nu$ is the bundle projection of the quaternionic Hopf bundle $Sp(1)\hookrightarrow S^7\!=\!\partial D^8\to\mathbb{H}P^1\!=\!S^4$.  
Then, observing this cellular decomposition of 
$\mathbb{H}P^2$ in view of the proof of (ii), we find that for $l\in\mathbb{Z}_+$, the null homotopy of $3l\,\tau_{Sp(3)/Sp(1), \,(\mathscr{L}^\lambda)_{Sp(1)}}$ 
can be deformed into a null homotopy of the stable map of $3l\,\nu$ for the reason that $[Sp(2), \mathscr{L}]$ is of order 3 since the fibration above has 
$Sp(2)$ as the fiber. This together with the fact that 
$\pi^S_{18}$ is a 2-group shows that 
$\tau_{Sp(3)/Sp(1), \,(\mathscr{L}^\lambda)_{Sp(1)}}$ must be of order $8$ since  
$\nu$ is of order $24$. Hence we have that it provides a generator of the first direct summand of $\pi^S_{18}=\mathbb{Z}_8a[18]\oplus\mathbb{Z}_2$; namely that we can take $[Sp(3)/Sp(1), (\mathscr{L}^\lambda)_{Sp(1)}]$ to be $a[18]$.  
\end{proof}

\begin{proof}[Proof of $\mathrm{(vii)}$]
Following the above case we consider again the fibration  
\[Sp(2)\hookrightarrow Sp(3)/Sp(1)_1\to Sp(3)/Sp(1)_1\times Sp(2)
=\mathbb{H}P^2\]
where $Sp(1)_1=Sp(1)$. Let us write $Sp(1)_2\times Sp(1)_3$ for $Sp(1)\times Sp(1)\subset Sp(2)$ and consider that the framing $\mathscr{L}^{2\lambda}$ on 
$Sp(3)$ behaves on each of $Sp(1)_1$ and $Sp(1)_2$ as $\mathscr{L}^\lambda$. Then by combining the argument in the case (ii) with that in the case (vi) we find that $\tau_{Sp(3)/Sp(1)_1, \,(\mathscr{L}^{2\lambda})_{Sp(1)_1}}$ must yields an 
element of order 2 in $\pi^S_{18}$ which is a 2-group due to the fact that  
$[Sp(2)/Sp(1)_1, (\mathscr{L}^\lambda)_{Sp(1)_1}]$ has order $240=10\cdot 24$. 
So we have 
\[\mathrm{ord}([Sp(3)/Sp(1)_1, (\mathscr{L}^{2\lambda})_{Sp(1)_1}])=2\]
where $\mathrm{ord}(a)$ denotes the order of $a$.

Let $\xi$ be the complex line bundle associated to the principal 
circle bundle    
\[S\hookrightarrow Sp(3)/Sp(1)_1\to Sp(3)/Sp(1)_1\times S\]
where $S\subset Sp(1)_2$ in the notation above
and let 
\[\pi\colon Sp(3)/Sp(1)_1\times S\to Sp(3)/Sp(1)_1\times Sp(2)=\mathbb{H}P^2\]
be the canonocal projection. Then observing $\xi$ under the decompositon 
$\mathbb{H}P^2=S^4\cup_\nu D^8$ we find that the restriction of it on 
$\pi^{-1}(\partial D^8)=S^7\times Sp(2)/S$ becomes stably trivial due to the action of $\lambda$ and therefore that $\xi$ itself is stably trivial. Based on this fact, applying the defining equation of Lemma 2.2 in ~\cite{O} to the sphere bundle $S(\xi)$ we have 
\[[Sp(3)/Sp(1), (\mathscr{L}^{2\lambda})_{Sp(1)}]=
\eta [Sp(3)/Sp(1)\times S, (\mathscr{L}^{2\lambda})_{Sp(1)\times S}].\]
This together with the above result on the non zeroness of this element  
tells us that we can take  
$a[17]=[Sp(3)/Sp(1)\times S, (\mathscr{L}^{2\lambda})_{Sp(1)\times S}]$. 
\end{proof}

\begin{proof}[Proof of $\mathrm{(viii)}$] 
Consider the canonical fibration
\[Sp(1)/S\times Sp(2)\hookrightarrow Sp(3)/S\to Sp(3)/Sp(1)\times Sp(2)=\mathbb{H}P^2.\] Since $\pi^S_{12}=0$, $[Sp(1)/S\times Sp(2), \mathscr{L}]=0$. Using this, an argument similar to that in the proof of (vi) above leads us immediately to the conclusion that  
$\tau_{Sp(3)/S, \\,\mathscr{L}}$ and $\nu$ have the same order 24. Therefore  
we have that $\tau_{Sp(3)/S, \\,\mathscr{L}}$ yields a generator of  
$\pi^S_{20}=\mathbb{Z}_{24}a[20]$, which at once implies that we can take 
$[Sp(3)/S, \mathscr{L}]$ to be $a[20]$.
\end{proof}

\begin{proof}[Proof of $\mathrm{(ix)}$] 
Consider the canonical fibration
\[Sp(1)/S_1\times Sp(2)/S_3\hookrightarrow Sp(3)/S_1\times S_3\to 
Sp(3)/Sp(1)\times Sp(2)=\mathbb{H}P^2\]
where we write $S_1$ for $S\subset Sp(1)$ and  
$Sp(1)_2\times Sp(1)_3$ for $Sp(1)\times Sp(1)\subset Sp(2)$. Let us view the framing $\mathscr{L}^\lambda$ of $Sp(3)$ as behaving on $Sp(1)_1$ and $Sp(1)_3$ as $\mathscr{L}$ and on $Sp(1)$ as $\mathscr{L}^\lambda$. Then applying a simiar argument to that in the case (vii) we see that 
$\tau_{Sp(3)/S_1\times S_3, (\mathscr{L}^\lambda)_{S_1\times S_3}}$ is not null 
homotopic and but twice that is null homotopic since  
$\mathrm{ord}([Sp(2)/Sp(1), \mathscr{L}^\lambda])=5\cdot 2\cdot 24$ and 
$11\cdot 24\,\pi^S_{19}=0$. So we have  
\[\mathrm{ord}([Sp(3)/S\times S, (\mathscr{L}^\lambda)_{S\times S}])=2.\]
Furthermore we find that it is impossible to decompose $\tau_{Sp(3)/S_1\times S_3, (\mathscr{L}^\lambda)_{S_1\times S_3}}$ into a sum of multiple maps. These imply that we can take 
$b[19]=[Sp(3)/S\times S, (\mathscr{L}^\lambda)_{S\times S}]$.
\end{proof}

\begin{proof}[Proof of $\mathrm{(x)}$]
Let us consider the canonical fibration  
\begin{equation*} 
Sp(1)\hookrightarrow Sp(3)/T^2 \to Sp(3)/Sp(1)\times T^2  
\end{equation*}
where $T^2=S\times S$. Then applying the converse of the procedure in 
the proof of (ii)  we find that not only 
$\tau_{Sp(3)/Sp(1)\times T^2, \,(\mathscr{L}^\lambda)_{Sp(1)\times T^2}}$ can be deformed into $\tau_{Sp(3)/T^2, \,(\mathscr{L}^\lambda)_{T^2}}$ but also a null homotopy of  $\tau_{Sp(3)/Sp(1)\times T^2, 
\,(\mathscr{L}^\lambda)_{Sp(1)\times T^2}}$ can be deformed into that of 
$\tau_{Sp(3)/T^2, \,(\mathscr{L}^\lambda)_{T^2}}$. Therefore from the result in the case (ix) above we see that 
\[ [Sp(3)/Sp(1)\times T^2, \,(\mathscr{L}^\lambda)_{Sp(1)\times T^2}]\ne 0.\]

Now since $Sp(1)\times T^2\subset Sp(3)$ is a subgroup of maximal rank we know that it is impossible for the defining equation of Lemma 2.2 in ~\cite{O} to be applied to $Sp(3)/Sp(1)\times T^2$, i.e. that 
$[Sp(3)/Sp(1)\times T^2, (\mathscr{L}^\lambda)_{Sp(1)\times T^2}]$ cannot be written in the form $\eta\alpha$ $(\alpha\in \pi^S_{15})$. This tells us that we are allowed to take $a[16]=[Sp(3)/Sp(1)\times T^2, (\mathscr{L}^\lambda)_{Sp(1)\times T^2}]$. 
\end{proof}

\end{document}